\documentclass[titlepage,11pt]{article}
\usepackage{latexsym}
\usepackage{amsfonts}
\usepackage{amsmath}
\usepackage{mathtools}
\usepackage{tikz}
\usepackage{comment}
\newcommand\blackslug{\hbox{\hskip 1pt \vrule width 4pt height 8pt depth 1.5pt
        \hskip 1pt}}
\newcommand\bbox{\hfill \quad \blackslug \bigbreak}

\def\CC{\hbox{-}\cdots\hbox{-}}
\def\LL{,\ldots,}

\title{The vertex sets of subtrees of a tree}
\author{
Maria Chudnovsky\thanks{Supported by NSF grant DMS-2348219 and AFOSR grant FA9550-22-1-0083.}\\
Princeton University,\\ Princeton, NJ 08544, USA
\and
Tung Nguyen\thanks{Supported by a Porter Ogden Jacobus Fellowship, and by AFOSR grant
FA9550-22-1-0234 and NSF grant  DMS-2154169.}\\
Princeton University,\\ Princeton, NJ 08544, USA
\and
Alex Scott\thanks{Supported by EPSRC grant EP/X013642/1}\\
University of Oxford, \\
Oxford, UK
\and
Paul Seymour\thanks{Supported by AFOSR grant
FA9550-22-1-0234, and by NSF grant DMS-2154169.}\\
Princeton University,\\ Princeton, NJ 08544, USA}

\date{January 19, 2025; revised \today}

\newtheorem{thm}{}[section]

\newcommand{\Proof}{\noindent{\bf Proof.}\ \ }

\begin{document}
\maketitle
\begin{abstract}
Let $\mathcal{F}$ be a set of subsets of a set $W$. When is there a tree $T$ with vertex set $W$ such that each member of $\mathcal{F}$
is the set of vertices of a subtree of $T$? It is necessary that 
$\mathcal{F}$ has the Helly property and the intersection
graph of $\mathcal{F}$ is chordal. We will show that these two necessary conditions are together sufficient in the finite case, and more generally, they are
sufficient if no element of $W$ belongs to infinitely many infinite sets in $\mathcal{F}$. 
\end{abstract}

\section{Introduction}
Graphs in this paper have no loops or parallel edgesi, and may be infinite. A graph is {\em chordal} if all its induced cycles have length three.
L. Sur\'anyi (see page 584 of~\cite{gyarfas}) and others~\cite{buneman, gavril, walter} showed that a finite graph is chordal if 
and only if
it is the intersection graph of a set of subtrees of a tree. (This is not always true for infinite graphs, as Halin~\cite{halin} showed.)

Here is a similar but different question:  let $\mathcal{F}$ be a set of subsets of a set $W$. When 
is there a tree $T$ with vertex set $W$ such that each member of $\mathcal{F}$
is the set of vertices of a subtree of $T$? There are two natural necessary conditions:
\begin{itemize}
\item {\bf (The chordal property)} The intersection graph of $\mathcal{F}$ is chordal.
\item {\bf (The finite Helly property)} For all $S_1\LL S_k\in \mathcal{F}$ with $k$ finite, if $S_i\cap S_j\ne \emptyset$ for
$1\le i< j\le k$, then $S_1\cap\cdots\cap S_k\ne \emptyset$.
\end{itemize}
We will prove that:
\begin{thm}\label{mainthm0}
Let $\mathcal{F}$ be a set of subsets of a set $W$. Suppose that $W$ is finite, or more generally, 
no element of $W$ belongs to infinitely many infinite members of $\mathcal{F}$. Then there a tree $T$ with vertex set $W$ such that each member of $\mathcal{F}$
is the set of vertices of a subtree of $T$, if and only if $\mathcal{F}$ has the chordal property and the finite Helly property.
\end{thm}
(This will follow from a more general result, \ref{strongthm} below.)
Let us see first that the two necessary conditions are not always sufficient:
\begin{thm}\label{counterex}
There is a set $\mathcal{F}$ of subsets of a set $W$, such that $\mathcal{F}$ has the chordal property and the finite Helly property,
and yet there is no tree $T$ with vertex set $W$ such that each member of $\mathcal{F}$ is the vertex set of a subtree of $T$.
\end{thm}
\Proof
Let $W$ be the set of non-negative integers, and let 
$$\mathcal{F}=\{\{i,i+1\}:i\ge 1\}\cup \{\{0,i,i+1,i+2...\}:i\ge 1\}.$$ 
This satisfies the two necessary conditions, but there is no corresponding tree. To see the latter, suppose that $T$ is a tree 
with $V(T)=W$ and all
members of $\mathcal{F}$ are vertex sets of subtrees. In particular, $(i,i+1)$ is an edge of $T$ for each $i\ge 1$, so there is only one
more edge in $T$, and it is incident with $0$; and for each $i\ge 1$, making $(0,i)$ an edge does not work, because
$\{0,i+1,i+2...\}$ is supposed to be the vertex set of a subtree. This proves \ref{counterex}.~\bbox

The set of \ref{counterex} was derived from an example of Halin~\cite{halin}. He gave a 
chordal graph $G$ that is not expressible as the intersection graph of a set of 
subtrees of a tree. Let $G$ be any such graph, let $W$ be the set of all maximal cliques of $G$ (a {\em clique} of a graph is 
a set of
pairwise adjacent vertices, not necessarily maximal), and for each $v\in V(G)$, let $S_v$ 
be the set of all members of $W$
that contain $v$. Then $\mathcal{F}=\{S_v:v\in V(G)\}$ satisfies \ref{counterex}, as can easily be checked. 

A strengthening of \ref{mainthm0} is proved in the next section. Later in the paper, we consider the same question with paths instead 
of trees. Tucker~\cite{tucker} solved this question for finite graphs, and we extend his result to infinite graphs. Then we apply
this to develop what we think is the proper notion of ``path-width'' for infinite graphs.

\section{The main result}

Let us say that $\mathcal{F}$ is {\em well-founded} if there is no countable sequence $(S_i:i\ge 1)$ of members of $\mathcal{F}$ such that
$S_1\cap S_2\cap \cdots \ne \emptyset$ and $S_1\cap \cdots \cap S_i\not\subseteq S_{i+1}$ for each $i\ge 1$.
(If $(S_i:i\ge 1)$ is such a sequence, then all the sets $S_i$ must be infinite: we will use this later.)
We observe that the example given above for \ref{counterex} is not well-founded.
Our main result is:

\begin{thm}\label{strongthm}
Let $\mathcal{F}$ be a set of subsets of a set $W$, such that $\mathcal{F}$ has the chordal property and the finite Helly property,
and such that $\mathcal{F}$ is well-founded.
Then there is a tree $T$ with vertex set $W$ such that each member of $\mathcal{F}$ is the vertex set of a subtree of $T$.
\end{thm}
\Proof If $W\notin \mathcal{F}$, we could replace $\mathcal{F}$ by $\mathcal{F}\cup\{W\}$, and the 
hypotheses would still hold; so we may assume that $W\in \mathcal{F}$.
Let us say a {\em fleet} is a set $\mathcal{F}$ of subsets of a set $W$, such that 
\begin{itemize}
\item $\mathcal{F}$ has the chordal property and the finite Helly property;
\item  $\mathcal{F}$ is well-founded; and
\item $W\in \mathcal{F}$.
\end{itemize}
We call the members of a fleet $\mathcal{F}$ the {\em ships} of $\mathcal{F}$.
Note that every ship $S$ is nonempty, by the finite Helly property applied to $\{S\}$.
A ship of cardinality two is called an {\em edge-ship}.
A fleet $\mathcal{F}'$ is an
{\em extension} 
of a fleet $\mathcal{F}$ if $\mathcal{F}\subseteq \mathcal{F}'$, and $\mathcal{F}'\setminus\mathcal{F}$ contains only 
edge-ships. In that case, if $\mathcal{F}$ is well-founded then so is $\mathcal{F}'$.
A fleet $\mathcal{F}$ is {\em maximal} if no extension of $\mathcal{F}$ is different from $\mathcal{F}$.
\\
\\
(1) {\em For every fleet $\mathcal{F}$, there is an extension of $\mathcal{F}$ that is maximal.}
\\
\\
Let $I$ be a set that is linearly ordered by some relation $<$, and for each $i\in I$ let $\mathcal{F}_i$ be a fleet, such that
for all distinct $i,j\in I$, if $i<j$ then $\mathcal{F}_j$ is an extension of $\mathcal{F}_i$. Let 
$\mathcal{F}=\bigcup_{i\in I}\mathcal{F}_i$. We claim that $\mathcal{F}$ is a fleet. For every finite set $\mathcal{R}$ 
of ships of $\mathcal{F}$ that pairwise intersect,  there exists $i\in I$ such that 
$\mathcal{R}\subseteq \mathcal{F}_i$ (since $\mathcal{R}$ is finite), and since $\mathcal{F}_i$ is a fleet, there is a vertex 
that belongs to every member of $\mathcal{R}$; and so $\mathcal{F}$ has the finite Helly property. Similarly it has the chordal 
property, and trivially all its ships of size different from two belong to each $\mathcal{F}_i$. Since $\mathcal{F}$ only differs from 
each $\mathcal{F}_i$ by the addition of edge-ships, it is well-founded. This proves that $\mathcal{F}$ is a fleet
extending each $\mathcal{F}_i\;(i\in I)$. From Zorn's lemma, this proves (1).

\bigskip

Consequently, to prove the theorem it suffices to prove it for maximal fleets, so we may assume that $\mathcal{F}$ is maximal.
A nonempty subset $X\subseteq W$ is {\em disconnected} if the set of edge-ships included in $X$ is the edge set of a disconnected graph
with vertex set $X$. 
\\
\\
(2) {\em If $W$ is not disconnected then the theorem holds.}
\\
\\
If $W$ is not disconnected, there is a tree $T$ with vertex set $W$ such that all its edges
are edge-ships of $\mathcal{F}$. Suppose that there is a ship $S$ that is not the vertex set of a tree in $T$.
Hence there are at least two components of $T[S]$; choose a minimal path $P$ of $T$ that joins two vertices in different components 
of $T[S]$. 
Thus $P$ has length at least two,
since the ends of $P$ are in different components of $T[S]$. Let $P$ have vertices $p_1\CC p_k$ in order; thus, $k\ge 3$,
and $p_1,p_k\in S$, and $p_2\LL p_{k-1}\notin S$. If $k=3$ then the three ships $\{p_1,p_2\},\{p_2,p_3\}, S$ pairwise intersect and yet
have no common vertex, contradicting that $\mathcal{F}$ has the finite Helly property; and if $k\ge 4$ then the intersection
graph of the set of ships $E(P)\cup \{S\}$ is a cycle of length at least four, contradicting the chordal property.
Thus 
every ship $S$ induces a tree in $T$, and so the theorem holds.
This proves (2).

\bigskip

Let us say a nonempty subset of $W$ expressible as an intersection of finitely many ships is a {\em meeting}. 
(Since $\mathcal{F}$ is well-founded, every nonempty intersection of ships is an intersection of finitely many ships, but we do not need that.)
If no meeting is disconnected, then the theorem holds by (2), since $W$ is a meeting;
so let us suppose that there is a disconnected meeting.
Since $\mathcal{F}$ is well-founded, there is no infinite sequence of meetings with nonempty intersection such that each is a proper
subset of its predecessor. 
It follows that there is a disconnected meeting $A$ such that no proper subset of $A$ is a disconnected meeting.
Let $F$ be the graph
with vertex set $A$ and edge set all edge-ships included in $A$. Thus $F$ is not connected; let the vertex sets of its components 
be $F_d\;(d\in D)$. 

Let $\mathcal{R}$ be the set of all ships that do not include $A$ as a subset. 
\\
\\
(3) {\em There is no sequence of ships $R_1\LL R_k\in \mathcal{R}$ such that for some distinct $d,d'\in D$, 
$R_1\cap F_d\ne \emptyset$  and $R_k\cap F_{d'}\ne \emptyset$, and 
$R_i\cap R_{i+1}\ne \emptyset$ for $1\le i<k$.}
\\
\\
Suppose there is such a sequence $R_1\LL R_k$,  and choose one with $k$ minimum. 
Choose finitely many ships $S_1\LL S_\ell$ such that $S_1\cap \cdots\cap S_{\ell} = A$.
If $k=1$, then $R_1\cap A$ is a disconnected 
meeting, and is a proper subset of $A$, a contradiction.
If $k=2$, then from the finite Helly property, since every two of the ships $S_1\LL S_\ell, R_1,R_2$ have nonempty intersection, they have a common 
member $c$ say. So $c\in A\cap R_1\cap R_2$, and so one of $R_1,R_2$ meets two of the sets $F_d\;(d\in D)$, contrary to the minimality of $k$.
Thus $k\ge 3$, and $R_i\cap R_j=\emptyset$ for $1\le i,j\le k$ with $j>i+1$; and $R_2\LL R_{k-1}$ are all disjoint from $A$ (because otherwise
we could reduce $k$). If each of $S_1\LL S_\ell$ has nonempty intersection
with $R_2$, then $S_1\cap\cdots \cap S_\ell\cap R_2$ is nonempty by the finite Helly property, contradicting that 
 $A\cap R_2= \emptyset$. So we assume that $S_1\cap R_2=\emptyset$. Since $S_1\cap R_k\supseteq A\cap R_k\ne \emptyset$, 
we may choose $j\in \{3\LL k\}$
minimum such that $S_1\cap R_j\ne \emptyset$. Thus $S_1$ has nonempty intersection with $R_1, R_j$ and is disjoint from $R_2\LL R_{j-1}$. 
Since $j\ge 3$, the intersection graph of the set of ships $\{R_1\LL R_j, S_1\}$ is a cycle of length at least four,
contradicting that $\mathcal{F}$ has the chordal property. This proves (3).

\bigskip

Choose distinct $d,d'\in D$ and $a\in F_d$ and $b\in F_{d'}$. Since $\{a,b\}$ is not an edge-ship, the maximality
of $\mathcal{F}$ tells us that $\mathcal{F}\cup \{\{a,b\}\}$ is not a fleet, and therefore violates either the chordal property
or the finite Helly property. 

Suppose that $\mathcal{F}\cup \{\{a,b\}\}$ does not have the finite Helly property. Thus 
there are finitely many ships $P_1\LL P_m$, pairwise with nonempty
intersection, and each containing one or both of $a,b$, such that neither of $a,b$ belongs to all of $P_1\LL P_m$. We may assume
that $a\notin P_1$ and $b\notin P_2$, and so $P_1,P_2\in \mathcal{R}$. Since $P_1\cap P_2\ne \emptyset$, this contradicts (3). 

Thus $\mathcal{F}\cup \{\{a,b\}\}$ does not have the chordal property; and so there are finitely many ships $P_1\LL P_m$
such that the intersection graph of $\{P_1\LL P_m,\{a,b\}\}$ is a cycle of length at least four. Consequently each of $P_1\LL P_m$
contains at most one of $a,b$, and so $P_1\LL P_m\in \mathcal{R}$, 
contrary to (3). 
This contradiction show that there is no disconnected meeting, and so the theorem holds by (2). This proves \ref{strongthm}.~\bbox

We can obtain a slight strengthening of the theorem of Halin~\cite{halin}. He proved that if $G$ is a chordal graph with no infinite clique, then $G$ is the intersection graph of a set of subtrees of a tree. We will prove:
\begin{thm}\label{betterhalin}
Let $G$ be a chordal graph, such that no countable clique can be ordered as $\{v_i:i\ge 1\}$ in such a way that for each $i\ge 1$,
there is a vertex adjacent to all of $v_1\LL v_{i-1}$ and not to $v_i$.
Then $G$ is the intersection graph of a set of subtrees of a tree.
\end{thm}
\Proof
Let $W$ be the set of all maximal cliques of $G$, and for each $v\in V(G)$, let $S_v$ be the set of all maximal cliques that 
contain $v$. From the hypothesis, $\mathcal{F}=\{S_v:v\in V(G)\}$ satisfies the hypotheses of \ref{strongthm}, and so there is a tree
$T$ with vertex set $W$ such that each member of $\mathcal{F}$ is the vertex set of a subtree of $T$. But then $G$ is the intersection graph of this set of trees. This proves \ref{betterhalin}.~\bbox

\section{Interval ships}

There is a natural variation on the question we just discussed for subtrees: given a set of subsets 
of a set $W$, when is there a linear ordering of $W$ such that each of the sets is an interval under this ordering? 
When $W$ is finite, this is closely related to a theorem of Lekkerkerker and Boland~\cite{leck} (their characterization of 
interval graphs using ``asteroidal triples''), and was answered by Tucker~\cite{tucker}. 
If $\mathcal{F}$ is a set of subsets of a set $W$, we say $X\subseteq W$ is {\em $\mathcal{F}$-connected} if for every partition $(A,B)$ of $X$ into two nonempty sets $A,B$, there exists $S\in \mathcal{F}$ with $S\subseteq X$ such that 
$S\cap A, S\cap B\ne \emptyset$. Tucker~\cite{tucker} proved:
\begin{thm}\label{tuckerthm}
Let $\mathcal{F}$ be a set of subsets of a finite set $W$. Then there is a path $P$ with vertex set $W$ such that each member
of $\mathcal{F}$ is the vertex set of a subpath of $P$, if and only if there do not exist distinct $v_1,v_2,v_3\in W$, and $\mathcal{F}$-connected subsets
$X_1,X_2,X_3$ of $W$, such that $v_i\in X_j$ if and only if $i\ne j$ for all $i,j\in \{1,2,3\}$. 
\end{thm}
Tucker expressed his result in terms of ``matrices with the consecutive ones property'' and their associated bipartite graphs, but 
it is easy to derive one form from the other.

We would like to point out that this has an extension when $W$ is infinite. In fact, unlike \ref{strongthm}, the natural extension to 
infinite $W$ needs no extra hypothesis (except we have to replace ``path'' by ``linear ordering''). We will prove:
\begin{thm}\label{bettertuckerthm}
Let $\mathcal{F}$ be a set of subsets of a set $W$. Then there is a linear ordering $<$ of $W$ such that each 
member
of $\mathcal{F}$ is an interval under this ordering, if and only if there do not exist distinct $v_1,v_2,v_3\in W$, and subsets
$X_1,X_2,X_3$, each connected in $\mathcal{F}$, such that $v_i\in X_j$ if and only if $i\ne j$ for all $i,j\in \{1,2,3\}$.
\end{thm}
\Proof The ``only if'' part is easy, and we prove ``if''. 
If we have a linear ordering $<$ of a set $X$, and $a,b\in X$ are distinct, then the ordering of $X$ induces 
an ordering 
of $\{a,b\}$, say $a<b$. In this case we say the pair $(a,b)$ is {\em induced} by $<$.
We will add more conditions on the linear ordering we hope for, of the form that certain ordered pairs are induced by $<$.

Let $\mathcal{T}$ be a set of ordered pairs of $W$. If $X\subseteq W$, a linear ordering of $X$ is  
{\em $(\mathcal{F},\mathcal{T})$-acceptable} if 
\begin{itemize}
\item for each $S\in \mathcal{F}$, if $S\cap X\ne \emptyset$ then $S\cap X$ is 
an interval under the ordering; and
\item every pair $(a,b)$ in $\mathcal{T}$ with $a,b\in X$ is induced by the ordering of $X$.
\end{itemize}
We say that a subset $X$ of $W$ is {\em $(\mathcal{F},\mathcal{T})$-feasible} if there is a 
$(\mathcal{F},\mathcal{T})$-acceptable
linear ordering of $X$.

By Tucker's theorem, we may assume that every finite subset
of $W$ is $(\mathcal{F},\emptyset)$-feasible. 
By Zorn's lemma we may choose the set $\mathcal{T}$ 
maximal such that every finite subset of $W$ is $(\mathcal{F},\mathcal{T})$-feasible. 
We claim that:
\\
\\
(1) {\em  Let $a,b\in W$
be distinct. Then $\mathcal{T}$ contains exactly one of $(a,b), (b,a)$.}
\\
\\
Since $\{a,b\}$ is $(\mathcal{F},\mathcal{T})$-feasible, at most one of $(a,b), (b,a)$ is in $\mathcal{T}$; 
suppose that neither are.
By adding $(a,b)$ to $\mathcal{T}$, forming $\mathcal{T}_1$ say, the maximality of $\mathcal{T}$ tells us that there is a finite 
subset $X_1\subseteq W$
that is not $(\mathcal{F},\mathcal{T}_1)$-feasible; and so $a,b\in X_1$, and no $(\mathcal{F},\mathcal{T})$-acceptable 
linear ordering of $X_1$ induces $(a,b)$. Define $X_2$ similarly, adding $(b,a)$.
But $X_1\cup X_2$ is finite, and so is 
$(\mathcal{F},\mathcal{T})$-feasible, a contradiction. 
This proves (1).

\bigskip
Say $a<b$ if $a,b\in W$ are distinct and $(a,b)\in \mathcal{T}$.
Since every finite subset of $W$ is $(\mathcal{F},\mathcal{T})$-feasible, it follows that if $a,b,c\in W$ then not all three of 
$(a,b),(b,c), (c,a)\in \mathcal{T}$; and so, by (1),  $<$ is a linear ordering of $W$. We claim that every $S\in \mathcal{F}$
is an interval under this ordering. Suppose that $a,b,c\in W$ where $a<b<c$ and $a,c\in S$ and $b\notin S$; then  $\{a,b,c\}$
is not $(\mathcal{F},\mathcal{T})$-feasible, a contradiction. 
This proves \ref{bettertuckerthm}.~\bbox

There are many similar questions. For instance, given a set of subsets $\mathcal{F}$
of a set $W$, when is there a tree $T$ with vertex set $W$ such that each member of $\mathcal{F}$ is the vertex set of a path of $T$? 
Even in 
the finite case, this is open as far as we know. Assuming the members of $\mathcal{F}$ are all nonempty, 
there are three necessary conditions we see: the finite Helly property and the chordal property, as in \ref{strongthm}, and in 
addition, for each $S_0\in \mathcal{F}$, the set of all sets of the form 
$S\cap S_0$ (for $S\in \mathcal{F}$) must satisfy the condition of \ref{tuckerthm}. We do not know whether these three are sufficient.

Another interesting question is, given a set of subsets $\mathcal{F}$
of a set $W$, when is there a tree $T$ with {\em edge set} $W$ such that each member of $\mathcal{F}$ is the edge set of a path of $T$?
In the finite case, this can be answered, as follows. For each $S\in \mathcal{S}$, take a new element $v_S$.
Let $V$ be the set of all the new
elements $v_S\;(S\in \mathcal{F})$, and let $M$ be the binary matroid with element set $W\cup V$, in which $W$ is a base and each set
$S\cup \{v_S\}$ is a fundamental circuit with respect to this base. The tree we want exists if and only if $M$ is graphic;
that is (by Tutte's famous theorem~\cite{tutte}) if and only if $M$ has no minor isomorphic to the Fano matroid or its dual, and no minor isomorphic to the bond matroid
of $K_5$ or of $K_{3,3}$. These ``minor'' conditions can presumably be reformulated as explicit conditions on $\mathcal{F}$, 
though we have not worked out exactly what they say.

\section{Path-width and line-width}

There is an application of these results to a question about path-width. 
A {\em path-decomposition} of a finite graph $G$ is a sequence $(B_1\LL B_n)$ of subsets of $V(G)$, such that:
\begin{itemize}
\item for each vertex $v$ of $G$, there exists $i\in \{1\LL n\}$ with $v\in B_i$;
\item for each edge $uv$ of $G$,  there exists $i\in \{1\LL n\}$ with $u,v\in B_i$;
\item for each vertex $v$ of $G$ and all $i,i',i''\in \{1\LL n\}$ with $i\le i'\le i''$, if $v\in B_i\cap B_{i''}$ then $v\in B_{i'}$.
\end{itemize}
 We call the sets $B_i$ {\em bags}.
Its {\em width} is the maximum of $|B_i|-1$ over $i\in \{1\LL n\}$, and the {\em path-width} of a finite graph is the minimum width
of a path-decomposition. 

What happens for infinite graphs? One could try to extend the definition of path-decomposition, using one- or two-way
infinite sequences of bags instead of finite sequences, but it is not satisfactory. For instance, with path-width defined this way, graphs with uncountably many vertices and no edges do not have finite path-width, and the same holds for the disjoint union of countably many 
infinite paths, and the disjoint union of countably many infinite ``stars'' (trees with one vertex and countably many neighbours). 
So graphs in which every finite subgraph has path-width at most $k$ may have no finite path-width. This contrasts with tree-width 
(we omit the definition): Thomas~\cite{thomas} proved that if every finite subgraph of $G$ has tree-width at most $k$ then $G$ has 
tree-width at most $k$. 

Here is a better way to do it. Say a {\em line} is a set with a linear order $<$, and a {\em line-decomposition} of $G$ is a family 
$(B_i:i\in L)$ where $L$ is a line, and each $B_i$ is a subset of $V(G)$, satisfying:
\begin{itemize}
\item for each vertex $v$ of $G$, there exists $i\in L$ with $v\in B_i$;
\item for each edge $uv$ of $G$,  there exists $i\in L$ with $u,v\in B_i$;
\item for each vertex $v$ of $G$ and all $i,i',i''\in L$ with $i\le i'\le i''$, if $v\in B_i\cap B_{i''}$ then $v\in B_{i'}$.
\end{itemize}
Define the {\em line-width} of $G$ to be the maximum of $|B_i|-1$ over all $i\in L$, if this maximum exists.
This definition was used in~\cite{coarsepw}, to increase the generality of a result about graphs with bounded path-width.

For finite graphs, path-width and line-width are the same, but for infinite graphs they might not be.
We will prove the analogue to Thomas' theorem:
\begin{thm}\label{pathwidth}
For every integer $k\ge 0$, if every finite subgraph of $G$ has path-width at most $k$ then $G$ has line-width at most $k$.
\end{thm}
\Proof Let $G$ have the property that 
every finite subgraph of $G$ has path-width at most $k$, and choose $G$ maximal (under adding edges) with this property: this is 
possibly by Zorn's lemma. It follows that for all $u,v$ that are distinct and nonadjacent, adding the edge $uv$ would violate the 
property, so there is a finite subgraph $G_{uv}$ containing both $u,v$, such that if we add the edge $uv$ to $G_{uv}$, it no
longer has path-width at most $k$. So in every path-decomposition $(B_1\LL B_n)$ of $G_{uv}$, there is no $i\in \{1\LL n\}$ with $u,v\in B_i$. 

We claim that every finite subgraph $G'$ of $G$ admits a path-decomposition of width at most $k$ such that every bag is a clique of $G'$.
Let $H$
be the union of $G'$ and the subgraphs $G_{uv}$ for all nonadjacent $u,v\in V(G')$. Then $H$ is finite, and so admits a 
path-decomposition $(B_1\LL B_n)$ of width at most $k$. For each nonadjacent pair $u,v\in V(G')$, 
$$(B_1\cap V(G_{uv})\LL B_n\cap V(G_{uv}))$$ 
is a path-decomposition of $G_{uv}$ of width at most $k$; and so from the property 
of $G_{uv}$, none of $B_1\LL B_n$ include both of $u,v$.
Consequently 
$$(B_1\cap V(G')\LL B_n\cap V(G'))$$ 
is a path-decomposition of $G'$ of width at most $k$,  and each $B_i\cap V(G')$ is a clique, as claimed.

Each clique of $G$ has cardinality at most $k+1$.
Let $W$ be the set of all maximal cliques of $G$, and for each $v\in V(G)$, let $S_v$ be the set of cliques in $W$ that contain $v$. 
It suffices to show there is a linear order of $W$ such that $S_v$
is an interval for each $v\in V(G)$, because this would give us a set of bags in a linear order, making a line-decomposition. 
By \ref{bettertuckerthm} applied to $\{S_v:v\in V(G)\}$, this is true if and only if 
for every finite subset $U\subseteq W$ and every finite subset $X\subseteq V(G)$, there is a 
linear order of $U$, such that $S_v\cap U$ is an interval of this order for each $v\in X$. But for each such
$U$ and $X$, there is a finite subgraph $G'$ of $G$ such that each clique in $U$ is a clique of $G'$, and 
$X\subseteq V(G')$. Since $G'$ is finite, it admits a path-decomposition $(B_1\LL B_n)$
of width at most $k$ such that each $B_i$ is a clique of $G'$. Each member of $U$ is a subset of some $B_i$, and therefore equals 
some $B_i$
from its maximality; and so $(B_1\LL B_n)$ induces a linear order of $U$ with the desired property.
This proves \ref{pathwidth}.~\bbox

Thomas' (unpublished) theorem~\cite{thomas} itself can be proved the same way~\cite{thomassen}, using \ref{strongthm} (and the fact that all cliques of $G$ are finite, so the 
corresponding family $\mathcal{F}$ is well-founded) or Halin's theorem~\cite{halin} (that every chordal graph with no infinite clique is the intersection graph of subtrees of a tree),  in place of \ref{bettertuckerthm}.

\section*{Acknowledgement} 
Thanks to Reinhard Diestel for discussions on the background to this material. ADS would like to thank the Department
of Mathematics at the University of Arizona, where some of this work was done.

\end{document}